# Deleting edges from Ramsey minimal examples - a paper written leaving the "scaffolding in place"
## by Robert Cowen
## Queens College (CUNY)

1. <u>Introduction</u>. A famous remark attributed to Carl F. Gauss is "no self-respecting architect leaves the scaffolding in place after completing the building." This style of presenting mathematical research is surely elegant and economical; however it gives the reader almost no insight into how the results were actually obtained. In this note, I will describe in some detail how I was led to make a conjecture and the process that finally led me to a proof of the conjecture. The result is in an area of mathematics called Ramsey Theory. Ramsey theory is a flourishing area of graph theory with an enormous number of difficult open problems (see [2]). I am certainly not an expert in Ramsey Theory, never having done any research in this area before, and was not intending to do any original research when I began this project. Rather, I wanted to use some problems from Ramsey Theory to illustrate using Boolean computation for a paper that I was writing for the *Mathematica* Journal[1]. Most of my research involves logic, set theory and some graph theory (see my webpage[2]). The idea for the *Mathematica* paper was to take some combinatiorial problems, in this case, some standard problems from Ramsey Theory, translate them into *Mathematica*'s logical language and use its "industrial strength" Boolean computational capability to "solve" these problems. I felt that *Mathematica*'s Boolean capability was under-utilized and wanted to encourage others to make use of it.

2. <u>Ramsey Theory Preliminaries</u>. If $n$ is a positive integer, $K_n$, denotes the complete graph on $n$ vertices; that is, the graph with $n$ vertices that contains every possible edge between these vertices. If $s, t$ are positive integers, the Ramsey number, $r(s, t)$ is the smallest integer $p$, such that if the edges of $K_n$ are colored either red or blue, there must be either a red $K_s$ or a blue $K_t$. Interchanging colors implies, $r(s, t) = r(t, s)$. It is a well-known theorem, due to Frank Ramsey, that $p$ exists, for each $s$ and $t$ (see, for example, [5]). If $r(s, t) = p$, we shal call $K_p$ a *Ramsey Minimal Example*. Thus if $K_p$ is a Ramsey Minimal Example, there is a red/blue edge coloring of $K_{p-1}$ without a red $K_s$ or a blue $K_t$. Thus, it is natural to ask: how many edges, $e$, must be removed from $K_p$ before we can red/blue color the edges without getting a red $K_s$ or a blue $K_t$? Since $K_{p-1}$ has such a coloring and and $K_p$ can be obtained from $K_{p-1}$ by adding a new vertex and p - 1 edges, $1 \leq e \leq p - 1$. It is well known that r(3,3) = 6 (see [3],[4],[5]) and Figure 1 shows that e =1, when p = 6, since neither the red or blue subgraph contains a triangle.

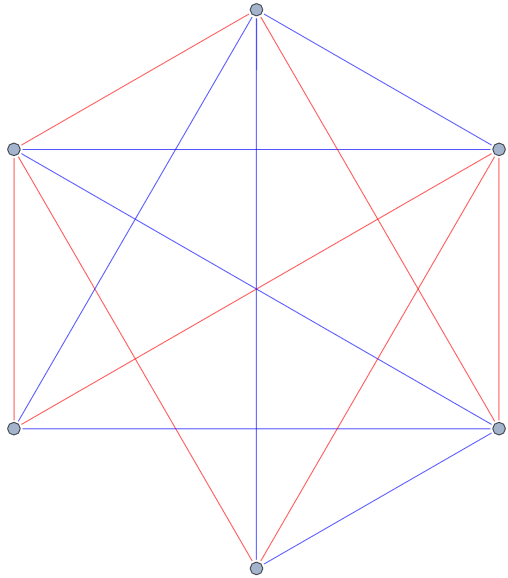

Figure 1

Utilizing Mathematica's Boolean capability I was able to show that in several other cases the answer also was $e=1$(some of these examples appear in [1]). I began to conjecture that this might always be true, but didn't have any ideas on how to prove or disprove it nor did I know if anyone had already solved this problem. I tried a "Google search" for this or related results. After several attempts, I entered the search terms: "removing edges ramsey theory" and the second entry displayed by the search was [3]. On this site I found a result from a paper by S. Golumb [4], attributed to his student, Herbert Taylor, that was a special case of the result I wanted to prove. The special case was when $K_s$ and $K_t$ are both triangles; however, the proof given in [4] was easily adapted to prove my conjecture.

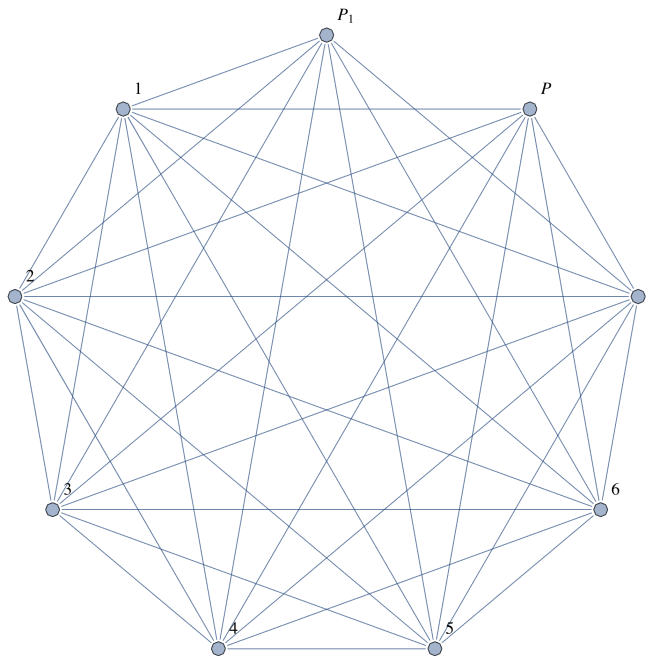

Figure 2

$K_p$

$K_s$    $K_t$

$K_{p-1}$    $P_1$    $P_1$    $K_{p-1}$

$P_1$    $L_{p-1}$    $K_{p-1}$    $P_1$    $P$

$P_1 Q$. Then $L_{p-1}$    $K_{p-1}$    $K_{p-1}$    $L_{p-1}$    $K_8$

Theorem. Suppose r(s, t) = $p$. If $G$ results from $K_p$ by deleting a single edge, then $G$ has a red/blue edge coloring and with respect to this coloring there is neither a red $K_s$ or a blue $K_t$.

Proof. Let $P$ be a vertex of $K_{p-1}$ and let $P_1$ be a new vertex. Connect $P_1$ to each vertex Q of $K_{p-1}$, except $P$, by a new edge, $P_1Q$. Let us denote by $L_{p-1}$, the graph that results from $K_{p-1}$ by replacing $P$ by $P_1$ and each edge $PQ$, $Q \neq P$, by $P_1 Q$. Then $L_{p-1}$ is isomorphic to $K_{p-1}$. (In Figure 2, $p = 9$, so $K_{p-1}$ and $L_{p-1}$ are isomorphic to $K_8$) Also, $G$ is isomorphic to $K_{p-1} \cup L_{p-1}$, the graph whose vertex set is the union of the vertex sets of $K_{p-1}$ and $L_{p-1}$ and whose edge set is the union of the edge sets of $K_{p-1}$ and $L_{p-1}$. We construct an edge coloring, $\sigma$, of $G$ as follows. Choose a red/blue edge coloring, $\gamma$ for $K_{p-1}$ that does not have either a red $K_s$ or a blue $K_t$. Extend this coloring to a coloring $\sigma$ of $K_{p-1} \cup L_{p-1}$, by coloring the additional edges, $P_1 Q$, the same color that $PQ$ received under $\gamma$. Then the subgraph $L_{p-1}$ also has the property that it does not contain a red $K_s$ or a blue $K_t$ with respect to the coloring $\sigma$. Any complete subgraph of $G$ cannot contain both $P$ and $P_1$, since there is no edge, $PP_1$ in $G$. Therefore, any $K_s$ that is a subgraph of $G$, is either a subgraph of $K_{p-1}$ (if it does not contain $P_1$) or a subgraph of $L_{p-1}$ (if it does not contain $P$). In either case $K_s$ is not colored red by $\sigma$. Similarly $K_t$ is not colored blue by $\sigma$.

Conclusion. I have found, in my own research, asking questions and doing computer "experiments," often leads to new theorems. Gauss, it should be remembered was a calculating prodigy and might have discovered many of his results this way as well. Unfortunately, by "removing the scaffolding," he also hid his methods of discovery.

<div align="center">References.</div>